\documentclass[a4paper, draft]{article}
\usepackage{amssymb,amsmath,mathrsfs,amsthm,amscd}

\theoremstyle{plain}
\newtheorem{theorem}{Theorem}[section]
\newtheorem*{conjecture}{Conjecture}

\newtheorem{lemma}[theorem]{Lemma}

\theoremstyle{definition}
\newtheorem{definition}[theorem]{ÎDefinition}

\theoremstyle{remark}
\newtheorem*{remark}{Remark}
\newtheorem*{example}{Example}
\newtheorem*{examples}{Examples}

\newcommand{\R}{\mathbb{R}}

\newcommand{\Q}{\mathbb{Q}}

\renewcommand{\le}{\leqslant}
\renewcommand{\ge}{\geqslant}

\newcommand{\Zs}{\mathcal{Z}}

\newcommand{\pd}{\partial}

\newcommand{\y}{{\bf y}}
\newcommand{\x}{{\bf x}}
\newcommand{\mdim}{\mathop{\rm mdim}}
\newcommand{\trk}{\mathop{\rm trk}}
\newcommand{\rk}{\mathop{\rm rk}}
\newcommand{\lk}{\mathop{\rm lk}}
\newcommand{\hrk}{\mathop{\rm hrk}}

\newcommand{\im}{\mathop{\rm im}}

\begin{document}
\title{Toral rank conjecture for moment-angle complexes}
\author{Ustinovsky Yury}
\date{}

\maketitle

\abstract{ We consider an operation $K\mapsto L(K)$ on the set of
simplicial complexes, which we call the ``doubling operation''. This
combinatorial operation has been recently brought into toric
topology by the work of Bahri, Bendersky, Cohen and Gitler on
generalised moment-angle complexes (also known as $K$-powers). The
crucial property of the doubling operation is that the moment-angle
complex $\mathcal Z_K$ can be identified with the real moment-angle
complex $\mathbb R\mathcal Z_{L(K)}$ for the double $L(K)$. As an
application we prove the toral rank conjecture for $\mathcal Z_K$ by
estimating the lower bound of the cohomology rank (with rational
coefficients) of real moment-angle complexes $\mathbb R\mathcal
Z_K$. This paper extends the results of our previous work, where the
doubling operation for polytopes was used to prove the toral rank
conjecture for moment-angle manifolds.\footnote{While preparing this
paper, the author learned that the toral rank conjecture for
moment-angle complexes has been also recently proved in a work of
Cao and Lu~\cite{ca-lu09}, albeit by a different method.} }

\section{Doubling operations}
Here we give the definition of the ``doubling operation'' and
discuss its main properties.

The author learned about the definition below from a communication
with the authors of the recent work~\cite{b-b-c-g}, in which the
importance of the ``doubling operation'' in toric topology has been
first demonstrated.
\begin{definition}
Let $K$ be an arbitrary simplicial complex on the vertex set
$[m]=\{v_1,\dots,v_m\}$. The \emph{double} of $K$ is the simplicial
complex $L(K)$ on the vertex set $[2m]=\{v_1, v_1',\dots,v_m,
v_m'\}$ determined by the following condition: $\omega\subset[2m]$
is the minimal (by inclusion) missing simplex of $L(K)$ iff $\omega$
is of the form $\{v_{i_1}, v_{i_1}',\dots, v_{i_k}, v_{i_k}'\}$,
where $\{v_{i_1},\dots, v_{i_k}\}$ is a missing simplex of $K$.
\end{definition}
If $K=\pd P^*$ is a boundary of the dual of the simple polytope $P$,
then $L(K)$ coincides with $L(P)^*$, see Definition 1 in \cite{U}.

The doubling operation for simple polytopes has also featured in the
recent work~\cite{gi-lo} of Gitler and Lopes de Medrano.
\begin{examples}
\
\begin{itemize}
\item If $K=\Delta^m$ is the $(m-1)$-dimensional simplex, then
$L(K)=\Delta^{2m}$.
\item If $K=\pd\Delta^m$ is the boundary of the $(m-1)$-dimensional
simplex, then $L(K)=\pd\Delta^{2m}$.
\end{itemize}
\end{examples}

It is easy to see that "doubling operation" respects join of the
simplicial complexes i.e. $L(K_1*K_2)=L(K_1)*L(K_2)$.

Given a simplicial complex $K$ we denote by $\mdim K $ the minimal
dimension of the maximal by inclusion simplices. Thus, for any $K$
$\mdim K\le\dim K$, and $K$ is pure iff $\mdim K=\dim K$.

The following lemma is the direct corollary from the definitions.
\begin{lemma}\label{mdim calc}
Let $K$ be a simplicial complex on $[m]$, then $\dim L(K)=m+\dim K$
and $\mdim L(K)=m+\mdim K.$
\end{lemma}
\section{$K$-powers}
\begin{definition}
Let $(X, A)$ be a pair of $CW$~--- complexes. For a subset
$\omega\subset[m]$ we define $$(X, A)^{\omega}:=\{(x_1,\dots,
x_m)\in X^m|x_i\in A \mbox{ for } i\not\in\omega\}.$$ Now let $K$ be
a simplicial complex on $[m]$. The $K$-\emph{power} of the pair $(X,
A)$ is $$(X, A)^K:=\bigcup_{\omega\in K}(X, A)^{\omega}.$$
\end{definition}
In this paper we shall consider two examples of $K$-powers (see
\cite{BPS}):
\begin{itemize}
\item Moment-angle complexes $\Zs_K=(D^2, S^1)^K.$
\item Real moment-angle complexes $\R\Zs_K=(I^1, S^0)^K.$
\end{itemize}

The next lemma explains the usefulness of the notion of "doubling
operation" in the studying of the relationship between moment-angle
complexes and real moment-angle complexes.
\begin{lemma}\label{Real - Ordinary lemma}
Let $(X, A)$ be a pair of $CW$~--- complexes and $K$ be a simplicial
complex on the vertex set $[m]$. Consider a pair $(Y, B)=(X\times X,
(X\times A)\cup (A\times X))$. For this pair we have:
$$(Y, B)^K=(X, A)^{L(K)}.$$ In particular $\Zs_K=\R\Zs_{L(K)}.$
\end{lemma}
\begin{proof}
For a point $\y=(y_1,\dots, y_m)\in Y^m$ we set
$$\omega_Y(\y)=\{v_i\in[m]\ |y_i\in Y\backslash B\}\subset[m].$$ For a point
$\x=(x_1, x_1',\dots, x_m, x_m')\in X^{2m}$ the subset
$\omega_X(\x)\subset[2m]$ is defined in a similar way. Let
$\y=(y_1,\dots, y_m)=((x_1, x_1'),\dots, (x_m, x_m'))\in
Y^m=X^{2m}$. It follows from the definition of the $K$-powers that
$\y\not\in (Y, B)^K$ iff $\omega_Y(\y)\not\in K$. The latter is
equivalent to the condition $\omega_X(\x)\not\in L(K)$, where
$\x=(x_1, x_1',\dots, x_m, x_m')$, since if
$\omega_Y(\y)=\{v_{i_1},\dots v_{i_k}\}$ then
$\omega_X(\x)\supset\{v_{i_1}, v_{i_1}',\dots, v_{i_k}, v_{i_k}'\}$.
Therefore $$\y\not\in (Y, B)^K \Leftrightarrow \x\not\in
(X,A)^{L(K)}$$ and the statement of the lemma is proved.
\end{proof}
\begin{example}
Let $K=\pd \Delta^2$ be the boundary of 1-simplex. Then we get
decomposition of 3-dimensional sphere: $$\Zs_K=D^2\times S^1\cup
S^1\times D^2=S^3.$$ On the other hand $L(K)=\pd \Delta^4$ and
$\R\Zs_{L(K)}=\pd I^4=S^3$ is the boundary of the standard
4-dimensional cube. So, in accordance with the lemma,
$\Zs_K=\R\Zs_{L(K)}$.
\end{example}
\section{Toral rank conjecture}
Let $X$ be a finite-dimensional topological space. Denote by
$\trk(X)$ the largest integer for which $X$ admits an \emph{almost
free} $T^{\trk(X)}$ action.
\begin{conjecture}[Halperin's toral rank conjecture, \cite{H}]
$$\hrk(X,\Q):=\sum \dim H^i(X, \Q)\ge 2^{\trk(X)}$$
\end{conjecture}
Moment-angle complexes provide a big class of spaces with torus
action, since there is natural coordinatewise $T^m$ action on the
space $\Zs_K$. In fact for some $r$ one can choose subtorus
$T^r\subset T^m$ such that the action $T^r\colon\Zs_K$ is almost
free. Our aim is to estimate the maximal rank of such subtorus and
the lower bound of $\hrk(\Zs_K, \Q)$.
\begin{lemma}
Let $K$ be $(n-1)$-dimensional simplicial complex on the vertex set
$[m]$. Then the rank of subtorus $T^r\subset T^m$ that acts almost
freely on $\Zs_K$ is less or equal to $m-n$.
\end{lemma}
\begin{proof}
For a subset $\omega\subset[m]$ we set $T^{\omega}=(T,e)^{\omega}$
(see definition of $K$-powers), where $e\in T$ is identity. It is
easy to see that isotropy subgroups of the action $T^m\colon \Zs_K$
are of the form $T^{\omega}$, $\omega\in K$. Therefore $T^r\subset
T^m$ acts almost freely iff the set $T^r \cap T^{\omega}$ is finite
for any $\omega\in K$.

Let $\sigma$ be the simplex of the dimension $(n-1)$. Since the
intersection $T^r \cap T^{\sigma}$ of two subtori in $T^m$ is
finite,
$$\rk T^r + \rk T^{\sigma}\le \rk T^m,$$ thus $r\le m-n$.
\end{proof}
\begin{remark}
In fact for any $(n-1)$-dimensional complex $K$ there is subtorus
$T^r\subset T^m$ of the rank $r=m-n$ that acts on $\Zs_K$ almost
freely, \cite[\S7.1]{DJ}.
\end{remark}
Now we prove our main result about the cohomology rank of the real
moment-angle complexes.
\begin{theorem}\label{hrk theorem}
Let $K$ be a simplicial complex on the vertex set $[m]$ with $\mdim
K=n-1$. Then $$\hrk(\R\Zs_K, \Q)\ge 2^{m-n}.$$
\end{theorem}
We first prove one general lemma.
\begin{lemma}\label{M-V lemma}

Let $(X,A)$ be a pair of $CW$-complexes such that $A$ has a collar neighbourhood
$U(A)$ in $X$, that is, $(U(A),A)\cong(A\times[0;1),A\times\{0\})$.
Let $Y=X_1\bigcup_A X_2$ be the space obtained by gluing two copies of $X$
along~$A$. Then the cohomology rank of $Y$ satisfies the inequality:
$$\hrk(Y,\Q)\ge \hrk(A,\Q).$$
\end{lemma}
\begin{proof}
Let $Y=X_1\cup X_2$ and $U_1(A)$, $U_2(A)$ be collar neighbourhoods
of $A$ in $X_1$ and $X_2$ respectively. Consider the covering
$Y=W_1\cup W_2$, where $W_1=X_1\cup U_2(A)$, $W_2=X_2\cup U_1(A)$
(thus $W_i\cong_{he}X$; $W_1\cap W_2=U_1(A)\cup U_2(A)\cong
A\times(-1,1)$). Now apply the Mayer Vietoris long exact sequence of
this covering to estimate the cohomology rank of $Y$:
$$
\begin{CD}
\cdots @>p_{(k-1)}^*>> H^{k-1}(W_1\cap W_2) @>\delta_{(k)}^*>>
H^{k}(Y) @>g_{(k)}^*>> H^{k}(W_1)\oplus H^k(W_2)
@>p_{(k)}^*>>\cdots
\end{CD}
$$
The map $p^*_{(k)}=(i_1^*)\oplus(-i_2^*)$ where $i_1$ and $i_2$ are
inclusions of $W_1\cap W_2$ in $W_1$ and $W_2$ respectively. Since
$W_1=W_2$ and the maps $i_1$ and $i_2$ coincide, $\dim
\ker{p^*_{(k)}}\ge\dim H^k(W_1)=\dim H^k(X)$. Applying this
inequality we have (remember that $W_1\cap W_2\cong_{he}A$):
\begin{align}\notag
&\dim H^k(Y)=\dim\ker g_{(k)}^*+\dim\im g_{(k)}^*=\dim\im
\delta_{(k)}^*+\dim\ker p_{(k)}^*\ge\dim H^{k-1}(A)-\\&-\dim\im
p_{(k-1)}^*+\dim H^k(X)\ge\dim H^{k-1}(A)-\dim H^{k-1}(X)+\dim H^k(X).\notag
\end{align}
After summing these inequalities over $k$ we obtain:
\begin{align}\notag
\hrk(Y,\Q)=\sum\dim H^k(Y)\ge\sum\bigl(\dim H^{k-1}(A)&-\dim H^{k-1}(X)+\dim
H^k(X)\bigr)=\\&=\sum\dim H^{k-1}(A)=\hrk(A,\Q).\notag
\end{align}

\end{proof}
\begin{proof}[Proof of the theorem 3.2]
We shall prove this fact by induction on $m$. The base of induction
is trivial.

Assume this statement is true for the complexes with less than $m$
vertices and $K$ is the complex with $m$ vertices.

The real moment-angle complex is a subspace of the $m$-dimensional
cube $\R\Zs_K\subset[-1;1]^m$. Denote by $(x_1,\dots, x_m)$
coordinates in $[-1;1]^m$. Assume that the vertex $v_1$ belong to
the maximal (by inclusion) simplex of $K$ of the dimension $\mdim
K=n-1$. Consider the decomposition of $\R\Zs_K=M_+\cup_X M_-$, where

$M_+=\{\vec{x}\in\R\Zs_K\subset\R^m\mid x_1\ge 0\}$,

$M_-=\{\vec{x}\in\R\Zs_K\subset\R^m\mid x_1\le 0\}$,

$X=\{\vec{x}\in\R\Zs_K\subset\R^m\mid x_1=0\}$.

It is easy to see that the pair $(M_+, X)$ satisfies the hypothesis
of the lemma \ref{M-V lemma}, so $$\hrk(\R\Zs_K,\Q)\ge \hrk(X,\Q).$$
Now lets describe the space $X$ more explicitly. Let $k$ be the
number of vertices in the complex $\lk v_1$. Then $X$ is just the
disjoint union of the $2^{m-k-1}$ copies of the space $\R\Zs_{\lk
v_1}$. Moreover, since $v_1$ is vertex of the maximal (by inclusion)
simplex of the minimal dimension $n-1$, so $\mdim\lk v_1=n-2$. Thus,
by the hypothesis of induction
$$\hrk(X,\Q)=2^{m-k-1}\hrk(\R\Zs_{\lk v_1},\Q)\ge 2^{m-k-1}\cdot2^{k-(n-1)}=2^{m-n}$$
The step of induction is proved.
\end{proof}
Now let's turn our attention to the moment-angle complexes.
Combining the results of lemma \ref{mdim calc}, lemma \ref{Real -
Ordinary lemma} and theorem \ref{hrk theorem} we have:
$$\hrk(\Zs_K,\Q)=\hrk(\R\Zs_{L(K)})\ge 2^{2m-\mdim L(K)-1}=2^{m-\mdim K-1}\ge2^{m-\dim K-1}.$$
Thus the \emph{toral rank conjecture} holds for the action of
subtori of $T^m$ on the moment-angle complexes $\Zs_K$.

The cohomology ring of $\Zs_K$ was calculated in \cite{BPS}. One of
the corollaries of this computation and Hochster's theorem states
(see \cite{BPS}, theorem 8.7):
\begin{theorem}\label{Z_K cohomology theorem}
$$H^*(\Zs_K,\Q)\cong \bigoplus_{\omega\subset[m],\
p\ge-1}\tilde{H}^p(K_{\omega},\Q),$$
where $K_\omega$ is the restriction of $K$ on the subset
$\omega\subset[m].$
\end{theorem}
In view of this theorem we can reformulate our main result as
follows:
$$\dim\bigoplus_{\omega\subset[m]}\tilde{H}^*(K_{\omega},\Q)\ge 2^{m-n},$$
for any simplicial complex $K$ on $[m]$ with $\mdim K=n-1$.

The author is grateful to V.\,M.~Buchstaber and scientific adviser
T.\,E.~Panov for suggesting the problem and attention to the
research.

\emph{E-mail address: }\texttt{yura.ust@gmail.com}
\end{document}